\documentclass[]{birkjour}
\usepackage[noadjust]{cite}
\usepackage{xcolor}
\RequirePackage[all]{xy}
\usepackage{graphicx}

\newtheorem{theorem}{Theorem}[section]
\newtheorem{lemma}{Lemma}[section]

\newtheorem{definition}{Definition}[section]

\newcommand{\BibTeX}{B\kern-0.1emi\kern-0.017emb\kern-0.15em\TeX}
\newcommand{\XYpic}{$\mathrm{X\kern-0.3em\raisebox{-0.18em}{Y}}$-$\mathrm{pic}\,$}

\newcommand{\cl}{C \kern -0.1em \ell}  

\newcommand{\ed}{\end{document}}
\begin{document}

%
%
%
%
%
%
%
%
%

\title[]{Quantum  solutions of a nonlinear Schr\"odinger equation}
\author[]{Sabrine Arfaoui}
%
\address{$^1$Laboratory of Algebra, Number Theory and Nonlinear Analysis, Department of Mathematics, Faculty of Sciences, University of Monastir,
Avenue of the Environment, 5019 Monastir, Tunisia.\\
$^2$Department of Mathematics, Faculty of Science, University of Tabuk, King Faisal Road, 47512 Tabuk, Saudi Arabia.}
\email{sabrine.arfaoui@issatm.rnu.tn}
%


\date{\today}
\dedicatory{Last Revised:\\ \today}
\begin{abstract}
In the present paper, we precisely conduct a $q$-calculus method for the numerical solutions of PDEs. A nonlinear Schr\"odinger equation is considered. Instead of the classical discretization methods we consider subdomains according to $q$-calculus, and provide an approximate solution due to a specific value of the parameter $q$. Error estimates show that $q$-calculus may produce efficient numerical solutions for PDEs.
\end{abstract}
\label{page:firstblob}
\maketitle
\section{Introduction}
The present paper is devoted essentially to the development of a numerical scheme to approximate the solution of a Nonlinear Schr\"odinger (NLS) equation in a quantum calculus framework. The aim crosses in fact the restriction to the resolution of such an equation, and goes further to show that $q$-calculus may provide good framework for numerical solutions of PDEs in general. It is well known that a major literature on the numerical solutions of PDEs is based on the classical methods such as finite difference, finite elements and volumes, Fourier analysis, and recently wavelets. See \cite{Benmabrouk-Ayadi1,Benmabrouk-Benmohamed-Omrani,Beziaetal,Bratsos1,Bratsos2,Bratsos-et-al1,Dehghan2006,Delfour-Fortin-Payre,Goncalves,Grillakis,Twizell-Bratsos-Newby}.

The NLS equation is in strong link to the modeling of real physical phenomena such as Newton's laws and conservation of energy in classical mechanics, behaviour of dynamical systems, the description of a particle in a non-relativistic setting in quantum mechanics, etc. Therefore, the NLS equation attracted the researchers from both theoretical and applied mathematics and physics. See \cite{Delfour-Fortin-Payre,Grillakis,Lamb,Twizell-Bratsos-Newby}.

Originally, Schr\"odinger's stated a linear form describing a moving particle according to the model equation
\begin{equation}\label{equa0}
\Delta\psi +\frac{8\pi^2m}{\hbar^2}\,\left(E-V(x)\right)\psi=0,
\end{equation}
where $\psi$ is known as a wave function, $m$ is the particle mass, $\hbar$ is the Planck's constant, $E$ is the energy, and $V$ is a potential. (\cite{byeon}, \cite{hase},\cite{malo},\cite{onor},\cite{zakh}).

Based upon the analogy between mechanics and optics, Schr\"odinger applied a perturbation method, to show an equivalence between his wave function in mechanics and Heisenberg's matrix. This gave rise next to the time dependent model
\begin{equation}\label{oheq}
i\hbar\psi_t=-\frac{\hbar^2}{2m}\,\Delta\psi+V(x)\psi-\gamma|\psi|^{2}\psi\qquad\mbox{in}\;\;\mathbb{R}^N,\;\;(N\geq2),
\end{equation}
known as the cubic NLS equation. 

Next, different variants and forms have been developed and investigated by researchers in different fields (\cite{abl},\cite{avron1},\cite{floer},\cite{kato},\cite{oh},\cite{sulem},). 

The present paper is devoted to the development of a numerical method based on $q$-calculus to approximate the solution of a reduced NLS equation in $\mathbb{R}^N$ written on the form
\begin{equation}\label{eqn1-1}
\left\{
\begin{array}{cc}
iu_{t}+\Delta\,u+f(u)=0,\quad\mbox{where}\;u=u(x,t)\;\;\mbox{with}\;(x,t)\in\Omega\times(t_0,+\infty)\\
u(x,t_0)=u_0(x)\quad\hbox{and}\quad\displaystyle\frac{\partial\,u}{\partial\,n}(x,t)=0,\quad(x,t)\in\Omega\times(t_0,+\infty).
\end{array}\right.
\end{equation}
We consider a domain $\Omega$ in $\mathbb{R}^N$ and $t_0$ a real parameter fixed as the initial time, $u_{t}$ is the first order partial derivative in time, $\Delta $ is the Laplace operator in $\mathbb{R}^N$. $\displaystyle\frac{\partial}{\partial n}$ is the outward normal derivative operator along the boundary $\partial\Omega$. $u$ and $u_0$ are complex valued functions. $f$ is a nonlinear function of $u$ assumed to at least continuous.

In \cite{chteouietal}, the stationary solutions of problem (\ref{eqn1-1}) has been studied using direct methods issued from the equation on the whole space. See also \cite{chteouietal1}. In \cite{Benmabrouk-Ayadi2} a Lyapunov-Sylvester method has been applied to solve numerical NLS and Heat equations. 
The organization of the present work will be as follows. In section 2, the $q$-calculus essential tools will be reviwed. Section 3 is devoted to the presentation of our main method. The discrete quantum version of a cubic NLS equation will be developed with necessary analysis of convergence, stability, solvability and consistency. Section 4 is subject of numerical experimentation due to our theoretical part. We conclude afterward.
\section{Quantum calculus toolkit}
One of the interesting fields of extensions of these analyses among other ones such as Hankel and Dunkel transforms is the so-called $q$-theory which is an important sub-field in harmonic analysis and which provides some discrete and/or some refinement of continuous harmonic analysis in sub-spaces such as $\mathbb{R}_q$ composed of the discrete grid $\pm q^n$, $n\in\mathbb{Z}$, $q\in(0,1)$. Recall that for all $x\in\mathbb{R}^*$ there exists a unique $n\in\mathbb{Z}$ such that $q^{n+1}<|x|\leq q^n$ which guarantees some density of the set $\mathbb{R}_q$ in $\mathbb{R}$. 

This section aims to introduce some basic concepts of $q$-theory. We present some definitions, notations and properties of $q$-derivatives and $q$-integrals which will be useful later. We are interested in some $q$-special functions which occupy a primordial place in this work namely the functions of $q$-Bessel ones. Backgrounds on $q$-theory may be found in \cite{Annaby-Mansour,Aral}, \cite{Dhaouadi},\cite{DhaouadiAtia},\cite{DhaouadiFitouhiElKamel},\cite{jackson,Koornwinder}  and the references therein.

For $0<q<1$, denote
$$
\mathbb{R}_{q}=\{\pm q^{n},\;\;n\in\mathbb{Z}\}\;\;\hbox{and}\;\;\widetilde{\mathbb{R}}_{q}^{+}=\mathbb{R}_{q}^{+}\bigcup \{0\}.
$$
We propose in this section to recall two basic functions that are applied almost everywhere in $q$-theory and its applications. See for example \cite{Gasper}.
\begin{definition}
The $q$-derivative of a function is defined by
$$
D_{q}f(x)=\begin{cases}
\dfrac{f(x)-f(qx)}{(1-q)x},\;\;&\;x\neq 0\\
f'(0)\;,& else,
\end{cases}
$$
provided that f is differentiable at 0.
\end{definition}
The operator $D_{q}f$ is the q-analogue of the classical derivative (\cite{Ali}). Indeed, if $f$ is differentiable, we get
$$
\lim_{q\longrightarrow 1}\,D_{q}f(x)=\frac{df(x)}{dx}.
$$
Many concepts of derivatives and integration rules  have been extended for the case of $q$-calculus. Many special functions have been also extended to especially Bessel, Exponential, Green, Mittag-Lefler functions. See \cite{Aral,jackson,Nadeem,Rezguietal}.

The only drawback of the $q$-calculus is the fact that they remain applied and investigated especially in harmonic functional analysis for the major part of the literature. A first step ahead has been conducted by Koornwinder and Swarttouw when studying Jackson's third $q$-Bessel function. Their work motivates researchers to develop different $q$-differential operators. Recently, $q$-calculus returns to take place in PDEs, indeed. Consider for example an elliptic equation
$$
\Delta u+f(u,x)=0
$$
where $f$ is a suitable function, generally nonlinear in $u$. We may search for a numerical $q$-approximation by considering a grid points in $\mathbb{R}_q$ instead of finite difference/finite elements used usually. In $q$-theory, we already have a $q$-analog of the Laplace operator expressed as
$$
\Delta_qu(x)=\dfrac{qu(q^{-1}x)-(1+q)u(x)+u(qx)}{x^2}.
$$
For $x=q^n$ in $\mathbb{R}_q^+$, we get
$$
u_{n+1}-(1+q)u_n+q u_{n-1}=-(1-q)^2qq^{2n}f_n
$$
where $u_n=u(q^n)$ and $f_n$ is some discretization of $f(u,x)$. We thus obtain a recursive equation permitting to compute $u_n$ recursively. More about applications of $q$-calculus in partial differential equations may be found in \cite{Annaby-Mansour}. A widely known example in $q$-theory is the Bessel type equation
$$
\begin{cases}
\Delta_{q}u(x)= -\lambda^{2}\,u(x),\\
u(0)=1,\;u'(0)=0,
\end{cases}
$$
($\lambda \in \mathbb{C}$), which has as unique solution a modified $q$-Bessel function. In the present paper, we will exploit the $q$-calculus to develop numerical solutions of some PDEs.
\section{The discrete NLS equation}
In this section we develop the details of our numerical method. 

For this aim we fix $\Omega=[0,1]$, $t_0=0$. Fix also a time step $l_k=(1-q)q^k$, and for $k\in\mathbb{N}$, we denote $t_k$ the $kth$ instant. For $n\in\mathbb{N}$, we denote $x_n=q^n$, and $h_n=q^n(1-q)$ the non-uniform space step. Denote also $u_n^k=u(t_k,x_n)$ the net function and $U_n^k$ its numerical approximation (the solution of the discrete problem). We discretize problem (...) as follows,
\begin{equation}\label{Discrete1}
i\dfrac{U_n^{k+1}-U_n^{k}}{l_k}+\dfrac{2q}{1+q}\dfrac{qU_{n-1}^k-(1+q)U_n^k+U_{n+1}^k}{h_{n}^2}+f(U_n^k)=0.
\end{equation}
By setting for $n,k\in\mathbb{N}$,
$$
\delta_n^k=\dfrac{2q}{1+q}\dfrac{l_k}{h_{n}^2},\qquad\hbox{and} \qquad\sigma_n^k=\dfrac{2q}{1+q}\dfrac{l_k}{h_{n}^2}i=\delta_n^ki, 
$$
the discrete problem (\ref{Discrete1}) becomes
\begin{equation}\label{Discrete2}
U_n^{k+1}=q\sigma_n^kU_{n-1}^k+(1-(1+q)\sigma_n^k)U_n^k+\sigma_n^kU_{n+1}^k+F_n^k,
\end{equation}
where $F_n^k=il_kf(U_n^k)$. Now denote for $k\in\mathbb{N}$,
$$
U^k=(U_n^k)_{n\in\mathbb{N}},
$$
the infinite vector of the numerical solution at the time $k$. Denote also
$$
\beta_n^k=1-(1+q)\sigma_n^k.
$$
We get the following dynamical infinite matrix-vector system
\begin{equation}\label{Discrete3}
U^{k+1}=A_kU^k+F^k,
\end{equation}
where $A_k$ is the infinite tri-diagonal matrix with coefficients
$$
A_k=\left(
\begin{array}{cccccccccc}
1-\sigma_0^k&\sigma_0^k&0&\dots&\dots&\dots&\dots&\dots&\dots\\
q\sigma_1^k&\beta_1^k&\sigma_1^k&0&\ddots&\ddots&\ddots&\ddots&\vdots\\
0&q\sigma_2^k&\beta_2^k&\sigma_2^k&0&\ddots&\ddots&\ddots&\vdots\\
0&0&q\sigma_3^k&\beta_3^k&\sigma_3^k&0&\ddots&\ddots&\vdots\\
\vdots&\ddots&\ddots&\ddots&\ddots&\ddots&\ddots&\ddots&\vdots\\
\vdots&\ddots&\ddots&\ddots&q\sigma_n^k&\beta_n^k&\sigma_n^k&\ddots&\vdots\\
\vdots&\ddots&\ddots&\ddots&\ddots&\ddots&\ddots&\ddots&\vdots\\
\vdots&\ddots&\ddots&\ddots&\ddots&\ddots&\ddots&\ddots&\vdots\\
\end{array}
\right)
$$
Infinite (especially tridiagonal) matrices are met in many fields and have been aplied widely. They are met in PDEs such as finite difference methods, in numerical analysis, and also in orthogonal polynomials theory and applications. These matrices appeared also in many physical problems such as optics and solid physics, quantum physics, etc. See for instance \cite{Akheizer}, \cite{Berezanskii}, \cite{Chihara}, \cite{Teschl}, \cite{Volkmer}.

In mathematics, infinite tridiagonal matrices are related to the so-called Jacobi operators (\cite{Teschl}). Spectral properties of these operators have been the subject of many studies, such as \cite{Dombrowski},\cite{Dombrowski1}, \cite{Dombrowski2} for the random case. An interesting problem in the theory of infinite matrices is the asymptotic behaviour of the eigen values spectrum and/or the characteristic polynomials. Such a question is widely met in ergodic theory for example. More about these matrices may be found in \cite{Damanik},\cite{Boutet},\cite{Janas},\cite{Janas1},\cite{Janas2},\cite{Janas3},\cite{Janas4},\\\cite{Janas5},\cite{Janas6},\cite{Janas7},\cite{Janas8},\cite{Janas9},\cite{Naboko},\cite{Naboko1},\cite{Naboko2}.

Notice in problem (\ref{Discrete3}) that a main difference with classical methods such as finite difference method is the possibility to relax one assumption on boundary conditions. We only need such an assumption for one extremity of the domain $\Omega$.

Now, observe that for each $k$, we get
$$
|1-\sigma_0^k|^2=1+|\delta_0^k|^2>|\delta_0^k|^2=|\sigma_0^k|^2,
$$
and similarly,
$$
|\beta_n^k|^2=1+(1+q)^2|\delta_n^k|^2>(1+q)^2|\delta_n^k|^2=(1+q)^2|\sigma_n^k|^2,
$$
which means that the matrix $A_k$ is a dominant-diagonal matrix, which guanrantees the solvability of our discrete scheme, and leads to the following theorem.
\begin{theorem}\label{thm-solvability}
The numerical problem (\ref{Discrete3}) is uniquely solvable, whenever $k\geq2n+1$.
\end{theorem}
In terms of the classical numerical schemes such as the finite difference, the assumption $k>>>n$ replaces the assumption $l=o(h^2)$, where $l$ and $h$ are the time and space steps for the finite difference scheme.

To investigate the stability of the numerical scheme, we propose to apply the Lyapunov criterion for stability, which states that a dynamical system $\mathcal{L}(U_{k},U_{k-1},\dots)=0$ is stable in the Lyapunov sense if for any bounded initial solution $U_{0}$ the solution $U_n$ remains bounded for all $n\geq0$ uniformly on $n$.
Here, we will precisely prove the following result.
\begin{lemma}\label{LyapunovStabilityLemma}
The solution $U^k$ is bounded independently of $k$ whenever the initial solution $U^0$ is bounded.
\end{lemma}
\textit{Proof.} 
We will proceed by recurrence on $k$. Assume firstly that $\|U^0\|\leq\eta$ for some $\eta$ positive. It follows from the fact that $k\geq2n+1$ that 
$$
|\sigma_n^k|\leq\dfrac{2q}{(1+q)(1-q)}.
$$
Therefore, using the system (\ref{Discrete2}), for $k=0$, we obtain
$$
U_n^{1}=q\sigma_n^0U_{n-1}^0+(1-(1+q)\sigma_n^0)U_n^0+\sigma_n^0U_{n+1}^0+F_n^0.
$$
As $U^0$ is bounded and the nonlinear function $f$ is continuous, we deduce that $U^1$ is bounded. So assume that $U^k$ is bounded. Using the system (\ref{Discrete2}) we get
$$
|U_n^{k+1}|\leq\dfrac{1+3q}{1-q}|U_n^k|+|F_n^k|.
$$
Using the recurrence hypothesis and again the continuity of the nonlinear function $F$, we deduce that
$$
|U_n^{k+1}|\leq\dfrac{2(1+q)}{1-q}C,
$$
where $C>0$ is a constant independent of $k$.

The consistency of the proposed method is done by evaluating the local truncation error arising from the discretization scheme. 
Assuming that the solution $u$ is sufficiently regular, we get the principal part as
\begin{equation}\label{consistency1}
\mathcal{L}(u)(x,t)=\dfrac{il_k}{2}u_{tt}+\dfrac{(1-q)h_n}{3q}u_{xxx}+o(l_k+h_n).
\end{equation}
It is clearly observable that the truncation operator $\mathcal{L}(u)$ goes to 0 as $n,k$ goes to infinity. This yields that the quantum numerical scheme is consistent at a minimum order 1 in time and space.

To finish with the convergence of the numerical method, we apply the Lax-Richtmyer equivalence theorem, which states that for consistent numerical approximations, stability and convergence are equivalent. We thus obtain the following lemma.
\begin{lemma}\label{laxequivresult}
As the numerical scheme is consistent and stable, it is then convergent.
\end{lemma}
Indeed, recall here that we have already proved in (\ref{consistency1}) that the used scheme is consistent. Next, Lemma \ref{LyapunovStabilityLemma}, yields the stability of the scheme. Consequently, the Lax equivalence Theorem guarantees the
convergence. So as Lemma \ref{laxequivresult}.
\section{Numerical implementations}
We propose in this experimental part to develop numerical examples to validate the theoretical results developed in the previous sections. We will use an $L_2$ discrete norm to evaluate the error between the exact solutions and the numerical ones as
\[
\|X\|_2=\Big(\sum_{i}|X_{i}|^2\Big)^{1/2},
\]
for any vector (series) $X=(X_{i})$ eventually in $L^2(\mathbb{C})$. Denote $u^k$ the net function $u(x,t^k)$ and $U^k$ the numerical solution. We propose to compute the discrete error
\begin{equation}\label{Er}
\mathrm{Er}=\max_k\|U^k-u^k\|_2
\end{equation}
on the grid $(x_{n})$, $n\geq0$.

We take for the rest $f(u)=|u|^2u$ which gives the original cubic NLS equation. We next take the classical soliton-type solution
$$
u(x,t)=\sqrt{\displaystyle\frac{2a}{q_s}}\exp\Bigl(i\bigl(\displaystyle\frac{1}{2}cx-\theta
t+\varphi\bigr)\Bigr)sech\Bigl(\sqrt{a}(x-ct)+\phi\Bigr)
$$
where $a$, $q_s$, $c$, $\theta=\displaystyle\frac{c^2}{4}-a$,
$\varphi$ and $\phi$ are some appropriate constants. For $t$
fixed, this function decays exponentially as
$|x|\rightarrow\infty$. It is a soliton-type disturbance which
travels with speed $c$ and with $a$-governed amplitude. See \cite{Benmabrouk-Ayadi1,Benmabrouk-Ayadi2,Bratsos1,Bratsos2,Bratsos-et-al1,hase,kato,Lamb,Twizell-Bratsos-Newby,zakh}.
\subsection{Propagation of a single soliton}
\begin{figure}[h]
\centering
\includegraphics[width=0.5\textwidth]{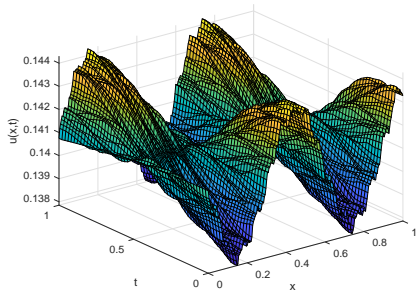}
\caption{Propagation of a single soliton, for $q=\frac{1}{8}$.}\label{Fig1SingleSoliton}
\end{figure}

\begin{figure}
\centering
\includegraphics[width=0.5\textwidth]{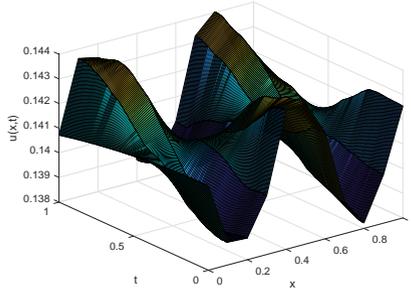}
\caption{Propagation of a single soliton, for $q=\frac{3}{8}$.}\label{Fig2SingleSoliton}
\end{figure}

In a first experimentation, we focus on a single-soliton-type particle. The computations are done for $0\leq x\leq 1$ and $0\leq t\leq 1$. We fix the $q$ parameter to many different values according to the closeness to 0 or to 1. Let $q\in\{q_i=\frac{i}{8},\;i=1,\dots,7\}$. We also fix the soliton parameters $a=0,01$, $q_s=1$, $c=0,1$ and the phase parameters $\varphi=\phi=0$. Figure \ref{Fig1SingleSoliton} and \ref{Fig2SingleSoliton} illustrate two cases of the numerical solution for the propagation of a single soliton issued from our quantum numerical scheme.
\newpage
Tables \ref{ErreurSingleSoliton50} and \ref{ErreurSingleSoliton100} illustrate the error estimates between the numerical solution and the exact one for different values of the quantum parameter $q$, and for different values of the maximal time index $K$. The space grid is fixed to a number $N$ of points. Truncating in an order $N$ for practice feasibility of the system (\ref{Discrete3}), we denote $TrU_N^k$ the truncated vector
$$
TrU_N^k=[U_0^k,U_1^k,\dots,U_N^k]^T,
$$
and $TrF_N^k$ the truncated vector
$$
TrF_N^k=[F_0^k,F_1^k,\dots,F_N^k]^T,
$$
where the upper script $^T$ is for the transpose. We denote similarly, $TrA^N_k$ the truncated matrix
$$
TrA_k^N=\left(
\begin{array}{cccccccccc}
1-\sigma_0^k&\sigma_0^k&0&\dots&\dots&\dots\\
q\sigma_1^k&\beta_1^k&\sigma_1^k&0&\ddots&\ddots\\
0&q\sigma_2^k&\beta_2^k&\sigma_2^k&0&\ddots\\
0&0&q\sigma_3^k&\beta_3^k&\sigma_3^k&0\\
\vdots&\ddots&\ddots&\ddots&\ddots&\ddots\\
\vdots&\ddots&\ddots&\ddots&q\sigma_N^k&\beta_N^k\\
\end{array}
\right).
$$
The system (\ref{Discrete3}) will be approximated by
$$
TrU_N^{k+1}=TrA^N_kTrU_N^k+TrF_N^k.
$$
\begin{table}
\centering
\begin{tabular}{llllllll}
\hline
\multicolumn{8}{c}{$K=10$}  \\
\hline
$q$ & $q_1$ & $q_2$ & $q_3$ & $q_4$ & $q_5$ & $q_6$ & $q_7$ \\
$Er$ & $1.17e-5$ & $1.28e-5$ & $2.41e-4$ & $2.52e-4$ & $2.18e-4$ & $3.11e-4$ & $2.87e-4$ \\
\hline
\multicolumn{8}{c}{$K=15$}  \\
\hline
$q$ & $q_1$ & $q_2$ & $q_3$ & $q_4$ & $q_5$ & $q_6$ & $q_7$ \\
$Er$ & $2.01e-6$ & $2.32e-6$ & $3.01e-5$ & $2.84e-5$ & $2.71e-5$ & $3.01e-5$ & $2.77e-5$ \\
\hline
\multicolumn{8}{c}{$K=20$}  \\
\hline
$q$ & $q_1$ & $q_2$ & $q_3$ & $q_4$ & $q_5$ & $q_6$ & $q_7$ \\
$Er$ & $1.27e-7$ & $1.44e-7$ & $1.23e-6$ & $2.15e-6$ & $2.53e-6$ & $2.01e-5$ & $2.76e-5$ \\
\hline
\end{tabular}
\caption{Error estimates for a single soliton for $N=20$.}
\label{ErreurSingleSoliton50}
\end{table}
\begin{table}
\centering
\begin{tabular}{llllllll}
\hline
\multicolumn{8}{c}{$K=10$}  \\
\hline
$q$ & $q_1$ & $q_2$ & $q_3$ & $q_4$ & $q_5$ & $q_6$ & $q_7$ \\
$Er$ & $1.17e-5$ & $1.28e-5$ & $2.41e-4$ & $2.52e-4$ & $2.18e-4$ & $3.11e-4$ & $2.87e-4$ \\
\hline
\multicolumn{8}{c}{$K=15$}  \\
\hline
$q$ & $q_1$ & $q_2$ & $q_3$ & $q_4$ & $q_5$ & $q_6$ & $q_7$ \\
$Er$ & $2.01e-6$ & $2.32e-6$ & $3.01e-5$ & $2.84e-5$ & $2.71e-5$ & $3.01e-5$ & $2.77e-5$ \\
\hline
\multicolumn{8}{c}{$K=20$}  \\
\hline
$q$ & $q_1$ & $q_2$ & $q_3$ & $q_4$ & $q_5$ & $q_6$ & $q_7$ \\
$Er$ & $1.27e-7$ & $1.44e-7$ & $1.23e-6$ & $2.15e-6$ & $2.53e-6$ & $2.01e-5$ & $2.76e-5$ \\
\hline
\end{tabular}
\caption{Error estimates for a single soliton for $N=50$.}
\label{ErreurSingleSoliton100}
\end{table}

It is noticeable from Table \ref{ErreurSingleSoliton50} and \ref{ErreurSingleSoliton100} that the numerical  quantum discrete scheme converges with good error estimates. This encorages to apply such a discretization idea for solving more complicated problems.
\subsection{Interaction of two solitons}
We consider here two solitons traveling with the same speed but in opposite directions in order to obtain the interaction phenomenon. The computations are done as in the previous experimentation on the space domain $\Omega_x=[0,1]$ and a time space also $I_t=[0,1]$. We fix the $q$ parameter as previously to $q\in\{q_i=\frac{i}{8},\;i=1,\dots,7\}$. We also fix the solitons parameters as follows,
\begin{itemize}
\item For the first soliton, we put $a=1$, $q_s=2$, $c=4$, $\varphi=0$ and $\phi=15$.
\item For the second soliton, we put $a=2,25$, $q_s=2$, $c=-4$, $\varphi=0$ and $\phi=-7,5$.
\end{itemize}
As in the previous experimentation, Tables \ref{ErreurTwoSolitons50} and \ref{ErreurTwoSolitons100} illustrate the error estimates between the numerical solution and the exact one for different values of the quantum parameter $q$, and for different values of the maximal time index $K$. The space grid is fixed to a number $N$ of points. 
\newpage
\begin{table}
\centering
\begin{tabular}{llllllll}
\hline
\multicolumn{8}{c}{$K=10$}  \\
\hline
$q$ & $q_1$ & $q_2$ & $q_3$ & $q_4$ & $q_5$ & $q_6$ & $q_7$ \\
$Er$ & $2.27e-5$ & $2.45e-5$ & $2.51e-4$ & $2.76e-4$ & $3.14e-4$ & $3.44e-4$ & $4.02e-4$ \\
\hline
\multicolumn{8}{c}{$K=15$}  \\
\hline
$q$ & $q_1$ & $q_2$ & $q_3$ & $q_4$ & $q_5$ & $q_6$ & $q_7$ \\
$Er$ & $2.11e-6$ & $2.18e-6$ & $2.05e-5$ & $2.61e-5$ & $2.85e-5$ & $3.12e-5$ & $4.05e-5$ \\
\hline
\multicolumn{8}{c}{$K=20$}  \\
\hline
$q$ & $q_1$ & $q_2$ & $q_3$ & $q_4$ & $q_5$ & $q_6$ & $q_7$ \\
$Er$ & $1.92e-7$ & $2.02e-7$ & $2.15e-6$ & $2.44e-6$ & $3.21e-6$ & $3.32e-5$ & $3.876e-5$ \\
\hline
\end{tabular}
\caption{Error estimates for two-interacted solitons for $N=20$.}\label{ErreurTwoSolitons50}
\end{table}
\begin{table}
\centering
\begin{tabular}{llllllll}
\hline
\multicolumn{8}{c}{$K=10$}  \\
\hline
$q$ & $q_1$ & $q_2$ & $q_3$ & $q_4$ & $q_5$ & $q_6$ & $q_7$ \\
$Er$ & $125e-5$ & $1.34e-5$ & $1.44e-4$ & $2.21e-4$ & $2.43e-4$ & $2.31e-4$ & $2.15e-4$ \\
\hline
\multicolumn{8}{c}{$K=15$}  \\
\hline
$q$ & $q_1$ & $q_2$ & $q_3$ & $q_4$ & $q_5$ & $q_6$ & $q_7$ \\
$Er$ & $1.15e-6$ & $1.27e-6$ & $1.32e-5$ & $1.17e-5$ & $1.95e-5$ & $2.02e-5$ & $2.33e-5$ \\
\hline
\multicolumn{8}{c}{$K=20$}  \\
\hline
$q$ & $q_1$ & $q_2$ & $q_3$ & $q_4$ & $q_5$ & $q_6$ & $q_7$ \\
$Er$ & $1.16e-7$ & $1.31e-7$ & $2.13e-6$ & $2.23e-6$ & $2.27e-6$ & $2.47e-5$ & $2.55e-5$ \\
\hline
\end{tabular}
\caption{Error estimates for two-interacted solitons for $N=50$.}\label{ErreurTwoSolitons100}
\end{table}

As in the previous case, we notice from Tables \ref{ErreurTwoSolitons50} and \ref{ErreurTwoSolitons100} that the numerical quantum discretization yielded a very close approximated solution to the exact one. This is clearly shown by the error estimates where the maximum error is estimated by $10^{-4}$ over all the values of the quantum parameter $q$ in the two tables \ref{ErreurTwoSolitons50} and \ref{ErreurTwoSolitons100}. This finding motivates the use of quantum numerical scheme for more general and/or complicated PDEs.
\section{Conclusion}
In the present paper, the principal aim was to test the efficiency of the quantum calculus in the approximation of the solutions of PDEs. As a prototypical example, we applied $q$-calculus to derive a numerical scheme for the well-known cubic NLS equation. As expected, the $q$-calculus yielded good approximations illustrated by low error estimates. The findings in the present paper make therefore good motivation to continue to exploit quamtum calculus for the numerical (and also exact) solutions of different types of PDEs. Comparisons with other models such as finite difference, finite volumes, and also wavelets as recent developments in mathematical analysis are fascinating and motivating future extensions. Copared to classical finite difference scheme method, we may conclude theoretically that the In fact, the present quantum scheme is more efficient, as it is based on geometric sequences time and space steps which surely converge rapidly than arithmetic discretizations. Therefore, we expect that involving or including hibrid schemes may induce best results. Finally, an interesting question rased from the present work may be formulated as follows: Given an infinite matrice that is truncated in an order $n$. We know that at most in $\mathbb{C}$, any truncation has at most $n$ eigenvalues. What can we expect for the original linear operator defined by means of the infinite matrice? This gives rise to possible chaotic behavior as a future study of the present case of matrices which are issued from parabolic, hyperbolic PDEs.. 


\begin{thebibliography}{10}
\bibitem{abl} M. J. Ablowitz, B.  Prinari, and A. D. Trubatch, Discrete and Continuous Nonlinear Schr\"odinger Systems. Cambridge Univ. Press, Cambridge, 2004.
	
\bibitem{Akheizer} N. I. Akheizer. The classical moment problem and some related questions in analysis. Oliver \& Boyd, 1965.
	
\bibitem{Annaby-Mansour} M. H. Annaby and Z. S. Mansour, $q$-Fractional Calculus and Equations. Lecture Notes in Mathematics 2056, Editors: J.-M. Morel and B. Teissier, Springer 2012.
	
\bibitem{Aral} A. Aral, V. Gupta and R. P. Agarwal, Applications of $q$-calculs in operator theory, Springer, New York, 2013.
	
\bibitem{avron1}  J. Avron, I. Herbst, and B.  Simon, Schr\"odinger operators with electromagnetic fields. III. Atoms in homogeneous magnetic field, {\it Commun. Math. Phys.} {\bf 79} (1981), 529-572.

\bibitem{Benmabrouk-Benmohamed-Omrani} A. Ben Mabrouk, M. L. Ben Mohamed and K. Omrani, Finite difference approximate solutions for a mixed sub-superlinear equation. J. Applied mathematics and computation 187 (2007), 1007-1016.

\bibitem{Benmabrouk-Ayadi1} A. Ben Mabrouk and M. Ayadi, A linearized finite-difference method for the solution of some mixed concave and convex nonlinear problems. Applied Mathematics and Computation 197 (2008), 1-10.

\bibitem{Benmabrouk-Ayadi2} A. Ben Mabrouk and M. Ayadi, Lyapunov type operators for numerical solutions of PDEs. Applied Mathematics and Computation 204 (2008), 395-407.

\bibitem{Berezanskii} Y. M. Berezanskii, Expansion in Eigenfunction of Self-Adjoint Operators, AMS, Providence, RI, 1968. Russian edition: Naukova Dumka, Kiev (1965).

\bibitem{Beziaetal} A. Bezia, A. Ben Mabrouk and K. Betina, Lyapunov-Sylvester Operators For $(2+1)$-Boussinesq Equation. Electronic Journal of Differential Equations, Vol. 2016 (2016), No. 268, pp. 1--19.

\bibitem{Bratsos1} A. G. Bratsos, A linearized finite-difference method for the solution of the nonlinear cubic Schr\"odinger equation, Comm. in Appl. Analysis 4(1) (2000), 133-139.

\bibitem{Bratsos2} A. G. Bratsos, A linearized finite-difference scheme for the numerical solution of the nonlinear cubic Schr\"odinger equation. Korean J. Comput. \& Appl. Math. 8(3) (2001), 459-467.

\bibitem{Bratsos-et-al1} A. G. Bratsos, Ch. Tsituras and D. G. Natsis, Linearized numerical schemes for the Boussinesq equation. Appl. NUm. Anal. Comp. Math. 2(1) (2005), 34-53.

\bibitem{byeon} J. Byeon and Z. Q. Wang, Standing waves with a critical frequency for nonlinear Schr\"odinger equations, {\it Arch. Ration. Mech. Analysis} {\bf 165} (2002), 295-316.

\bibitem{Chihara} T. S. Chihara, An introduction to orthogonal polynomials. In (1978).

\bibitem{chteouietal} R. Chteoui, A. Ben Mabrouk and H. Ounaiess, Existence and Properties of Radial Solutions of a Sub-linear Elliptic Equation. J. Part. Diff. Eq. 28 (1) (2015), 30-38

\bibitem{chteouietal1} R. Chteoui and A. Ben Mabrouk, A Generalized Lyapunov-Syslvester Computational Method for Numerical Solutions of NLS Equation With Singular Potential. Anal. Theory Appl., 33 (2017), pp. 333-354.

\bibitem{Damanik} D. Damanik and S. Naboko, Unbounded Jacobi matrices at critical coupling. Journal of approximation theory 145(2) (2007), p. 221-236.

\bibitem{Boutet} A. B. De Monvel, J. Janas and S. Naboko, Unbounded Jacobi matrices with a few gaps in the essential spectrum : constructive examples. Integral Equations and Operator Theory 69(2) (2011), p. 151-170.

\bibitem{Dhaouadi} L. Dhaouadi, On the $q$-Bessel Fourier transform, Bulletin of mathematical analysis and applications, 5(2) (2013), pp. 42-60.

\bibitem{DhaouadiAtia} L. Dhaouadi and M. J. Atia, Jacobi operators, $q$-dfference equations and orthogonal polynomials, arXiv:1211.0359v1, 2 Nov 2012, 22 pages.

\bibitem{DhaouadiFitouhiElKamel} L. Dhaouadi, A. Fitouhi and J. El Kamel, Inequalities in $q$-Fourier analysis, J. of Inequalities in Pure and Applied Mathematics, 7(5) (2006), Article 171, 14 pages.

\bibitem{Dehghan2006} M. Dehghan, Finite difference procedures for solving a problem arising in modeling and design of certain optoelectronic devices, Mathematics and Computers in Simulation 71 (2006), 16--30.

\bibitem{Delfour-Fortin-Payre} M. Delfour, M. Fortin and G. Payre, Finite difference solutions of a non-linear Schr\"odinger equation. J. Computa. Phys. 44 (1981), 277-288.

\bibitem{Dombrowski} J. Dombrowski and S. Pedersen, Spectral measures and Jacobi matrices related to Laguerre-type systems of orthogonal polynomials. Constructive approximation 13(3) (1997), p. 421-433.

\bibitem{Dombrowski1} J. Dombrowski, Tridiagonal matrix representations of cyclic self-adjoint operators II. Pacific Journal of Mathematics 120(1) (1985), p. 47-53.

\bibitem{Dombrowski2} J. Dombrowski, Tridiagonal matrix representations of cyclic selfadjoint operators. Pacific Journal of Mathematics 114(2) (1984), p. 325-334.

\bibitem{floer} A. Floer and A. Weinstein, Nonspreading wave packets for the cubic Schr\"odinger equation with a bounded potential, {\it J.~Funct. Anal.} {\bf 69} (1986), 397-408.

\bibitem{Gasper} G. Gasper and M. Rahman, Basic Hypergeometric serie Second edition, Combridge university Press, 2004.

\bibitem{Goncalves} E. Gon\c calv\`es, Resolution numerique, discretisation des EDP et EDO. Cours, Institut National Polytechnique de Grenoble, 2005.

\bibitem{Grillakis} M. G. Grillakis, On nonlinear Schr\"odinger equations. Commun. Partial. Differ. Equations. 25 (2000), 1827-1844.

\bibitem{hase} A. Hasegawa and Y. Kodama, Solitons in Optical Communications. Academic Press, San Diego, 1995.

\bibitem{jackson} F. H. Jackson, The application of basic numbers to Bessel's and Legendre's functions, Proc. London math.Soc. (2) 2 (1905) 192-220.

\bibitem{Janas} J. Janas and S. Naboko, On the point spectrum of some Jacobi matrices. Journal of Operator Theory (1998), p. 113-132.

\bibitem{Janas1} J. Janas and S. Naboko, Jacobi matrices with absolutely continuous spectrum. Proceedings of the American Mathematical Society 127(3) (1999), p. 791-800.

\bibitem{Janas2} J. Janas and S. Naboko, Jacobi matrices with power-like weights—grouping in blocks approach. Journal of Functional Analysis 166(2) (1999), p. 218-243.

\bibitem{Janas3} J. Janas and S. Naboko, Multithreshold spectral phase transition examples in a class of unbounded Jacobi matrices II. Citeseer, 2000.

\bibitem{Janas4} J. Janas and S. Naboko, Asymptotics of generalized eigenvectors for unbounded Jacobi matrices with power-like weights, Pauli matrices commutation relations and Cesaro averaging. Differential operators and related topics. Springer, 2000, p. 165-186.

\bibitem{Janas5} J. Janas and S. Naboko, Spectral analysis of selfadjoint Jacobi matrices with periodically modulated entries. Journal of Functional Analysis 191(2) (2002), p. 318-342.

\bibitem{Janas6} J. Janas and S. Naboko, Spectral properties of selfadjoint Jacobi matrices coming from birth and death processes. Recent Advances in Operator Theory and Related Topics. Springer, 2001, p. 387-397.

\bibitem{Janas7} J. Janas, S. Naboko and G. Stolz, Decay bounds on eigenfunctions and the singular spectrum of unbounded Jacobi matrices . In: International Mathematics Research Notices 4 (2009), p. 736-764.

\bibitem{Janas8} J. Janas, S. Naboko and G. Stolz, Spectral theory for a class of periodically perturbed unbounded Jacobi matrices : elementary methods. Journal of computational and applied mathematics 171(1-2) (2004), p. 265-276.

\bibitem{Janas9} J. Janas and S. Naboko, Criteria for semiboundedness in a class of unbounded Jacobi operators. Algebra i Analiz 14(3) (2002), p. 158-168.

\bibitem{kato} T. Kato, Remarks on holomorphic families of Schr\"odinger and Dirac operators, in {\it Differential equations}, (Knowles I., Lewis R., Eds.), North-Holland Math. Stud., Vol. 92, North-Holland, Amsterdam, 1984, pp. 341-352.

\bibitem{Koornwinder} T. H. Koornwinder and R. F. Swarttow, On $q$-Analogues of the Hankel and Fourier transform, Trans. A.M.S., 1992, 333, 445--461. 

\bibitem{Lamb} G. L. Lamb, Elements of soliton theory. Wiley 1980.

\bibitem{malo} B. A. Malomed, Variational methods in nonlinear fiber optics and related fields, {\it Progress in Optics} {\bf 43} (2002), 69-191.

\bibitem{Naboko} S. Naboko, I. Pchelintseva and L. O. Silva, Discrete spectrum in a critical coupling case of Jacobi matrices with spectral phase transitions by uniform asymptotic analysis. Journal of Approximation Theory 161(1) (2009), p. 314-336.

\bibitem{Naboko1} S. Naboko and S. Simonov, Spectral analysis of a class of hermitian Jacobi matrices in a critical (double root) hyperbolic case. Proceedings of the Edinburgh Mathematical Society 53(1) (2010), p. 239-254.

\bibitem{Naboko2} S. Naboko and S. Simonov. Titchmarsh-Weyl formula for the spectral density of a class of Jacobi matrices in the critical case. arXiv: 1911.10282 (2019)

\bibitem{Nadeem} R. Nadeem, T. Usman, K.S. Nisar, et al. A new generalization of Mittag-Leffler function via q-calculus. Adv Differ Equ 2020, 695 (2020). https://doi.org/10.1186/s13662-020-03157-z

\bibitem{oh} Y. G. Oh, Existence of semi-classical bound states of nonlinear Schr\"odinger equations with potentials of the class ($V_a$), {\it Communications in Partial Differential Equations} {\bf 13} (1988), 1499-1519.

\bibitem{onor} M. Onorato,A. R. Osborne, M. Serio, and S. Bertone, Freak waves in random oceanic sea states, {\it Phys. Rev. Lett.} {\bf 86} (2001), 5831-5834.

\bibitem{Rezguietal} I. Rezgui and A. Ben Mabrouk, Some Generalized $q$-Bessel type Wavelets and associated transforms, Anal. Theory Appl, 34(1) (2017), pp. 1-15.

\bibitem{sulem} C. Sulem and P.-L. Sulem, {The Nonlinear Schr\"odinger Equation. Self-focusing and Wave Collapse}, Applied Mathematical Sciences, Vol. 139, Springer-Verlag, New York, 1999.

\bibitem{Teschl} G. Teschl, Jacobi operators and completely integrable nonlinear lattices. T. 72. Mathematical Surveys and Monographs. American Mathematical Society, 2000, p. xvii+351.

\bibitem{Twizell-Bratsos-Newby} E. H. Twizell, A. G. Bratsos and J. C. Newby, A finite-difference method for solving the cubic Schr\"odinger equation. Mathematics \& Computers in Simulation. 43 (1997), 67-75.

\bibitem{Volkmer} H. Volkmer, Error estimates for Rayleigh-Ritz approximations of eigenvalues and eigenfunctions of the Mathieu and spheroidal wave equation. Constr. Approx. 20(1) (2004), p. 39-54.

\bibitem{zakh} V. E. Zakharov, {Collapse and Self-focusing of Langmuir Waves}, Handbook of Plasma Physics, (M. N. Rosenbluth and R. Z. Sagdeev, eds.), vol. 2 (A. A. Galeev and R. N. Sudan, eds.) 81-121, Elsevier (1984).

\end{thebibliography}
\end{document}